\documentclass[12pt]{amsart}

\newtheorem{theorem}{Theorem}[section]
\newtheorem{lemma}[theorem]{Lemma}
\newtheorem{proposition}[theorem]{Proposition}
\newtheorem{corollary}[theorem]{Corollary}

\theoremstyle{remark}
\newtheorem{definition}{Definition}
\newtheorem*{remark}{Remark}

\newcommand{\et}{\quad\mbox{and}\quad}
\newcommand{\bC}{\mathbb{C}}
\newcommand{\bN}{\mathbb{N}}

\newcommand{\bQ}{\mathbb{Q}}
\newcommand{\bR}{\mathbb{R}}
\newcommand{\bZ}{\mathbb{Z}}
\newcommand{\cA}{{\mathcal{A}}}

\newcommand{\cP}{{\mathcal{P}}}
\newcommand{\cO}{{\mathcal{O}}}
\newcommand{\cS}{{\mathcal{S}}}
\newcommand{\cU}{{\mathcal{U}}}
\newcommand{\ggll}{\gg\ll}
\newcommand{\GL}{\mathrm{GL}}
\newcommand{\MdR}{\mathrm{Mat}_{2\times2}(\bR)}
\newcommand{\MdZ}{\mathrm{Mat}_{2\times2}(\bZ)}
\newcommand{\red}{\mathrm{red}}

\newcommand{\uy}{\mathbf{y}}

\newcommand{\tA}{{\,{^t}\hskip-2pt A}}
\newcommand{\tB}{{\,{^t}\hskip-2pt B}}
\newcommand{\tL}{{\,{^t}\hskip-2pt L}}
\newcommand{\tM}{{\,{^t}\hskip-2pt M}}
\newcommand{\tN}{{\,{^t}\hskip-2pt N}}
\newcommand{\tU}{{\,{^t}\hskip-0.5pt U}}
\newcommand{\tW}{{\,{^t}\hskip-1pt W}}

\begin{document}

\baselineskip=17pt

\title[On the continued fraction expansion of a class of numbers]
{On the continued fraction expansion\\ of a class of numbers}

\author{Damien ROY}
\address{
   D\'epartement de Math\'ematiques\\
   Universit\'e d'Ottawa\\
   585 King Edward\\
   Ottawa, Ontario K1N 6N5, Canada}
\email{droy@uottawa.ca}
\dedicatory{ Au Professeur Wolfgang Schmidt,\\
avec mes meilleurs v{\oe}ux et toute mon admiration.}
\thanks{Work partially supported by NSERC and CICMA}

\maketitle


\section{Introduction}
\label{sec-intro}


A classical result of Dirichlet asserts that, for each real number
$\xi$ and each real $X\ge 1$, there exists a pair of integers
$(x_0,x_1)$ satisfying
\begin{equation*}
 1\le x_0\le X
 \et
 |x_0\xi-x_1|\le X^{-1}
\end{equation*}
(a general reference is Chapter I of \cite{Sc}). If $\xi$ is
irrational, then, by letting $X$ tend to infinity, this provides
infinitely many rational numbers $x_1/x_0$ with $|\xi-x_1/x_0|\le
x_0^{-2}$.  By contrast, an irrational real number $\xi$ is said to
be \emph{badly approximable} if there exists a constant $c_1>0$ such
that $|\xi-p/q| > c_1q^{-2}$ for each $p/q\in\bQ$ or, equivalently,
if $\xi$ has bounded partial quotients in its continued fraction
expansion.  Thanks to H.~Davenport and W.~M.~Schmidt, the badly
approximable real numbers can also be described as those
$\xi\in\bR\setminus\bQ$ for which the result of Dirichlet can be
improved in the sense that there exists a constant $c_2<1$ such that
the inequalities $1\le x_0\le X $ and $|x_0\xi-x_1|\le c_2X^{-1}$
admit a solution $(x_0,x_1)\in\bZ^2$ for each sufficiently large $X$
(see Theorem 1 of \cite{DSa}).

If $\xi$ is rational or quadratic real, then, upon writing
$\xi^2=(q\xi+r)/p$ for integers $p$, $q$ and $r$ with $p\neq 0$ and
putting $c_3=|p|\max\{|p|,|q|\}$, one deduces from the result of
Dirichlet that, for each $X\ge 1$, there exists a point
$(x_0,x_1,x_2) \in \bZ^3$ satisfying
\begin{equation*}
 1\le x_0 \le X,
 \quad
 |x_0\xi-x_1|\le c_3X^{-1}
 \et
 |x_0\xi^2-x_2|\le c_3X^{-1}.
\end{equation*}
Conversely, Davenport and Schmidt proved that, for each real number
$\xi$ which is neither rational nor quadratic over $\bQ$, there is a
constant $c_4>0$ such that, upon writing $\gamma = (1+\sqrt{5})/2$,
the system of inequations
\begin{equation}
 \label{intro:eq1}
 | x_0 | \le X, \quad
 | x_0\xi-x_1 | \le c_4X^{-1/\gamma}, \quad
 | x_0\xi^2-x_2 | \le c_4X^{-1/\gamma},
\end{equation}
admits no non-zero integer solution $(x_0,x_1,x_2)\in \bZ^3$ for
arbitrarily large values of $X$ (Theorem 1a of \cite{DS}). Since
$1/\gamma\simeq 0.618 <1$, this establishes a clear gap between the
set of rational or quadratic real numbers and the remaining real
numbers.  Moreover, this result of Davenport and Schmidt is best
possible in the following sense.  There exist real numbers $\xi$
which are neither rational nor quadratic and for which there is a
constant $c_5>0$ such that the system \eqref{intro:eq1}, with $c_4$
replaced by $c_5$, admits a non-zero integer solution for each $X\ge
1$ (Theorem 1.1 of \cite{Rb}).   These real numbers, which we call
\emph{extremal}, present from this point of view a closest behavior
to quadratic real numbers.  An application of Schmidt's subspace
theorem proves them to be transcendental over $\bQ$ (see Theorem 1B
in Chapter VI of \cite{Sc}).  Still they possess several properties
that make them resemble to quadratic real numbers. In the present
paper, we are interested in their approximation by rational numbers.

It is well known that each quadratic real number has an ultimately
periodic continued fraction expansion and so is badly approximable.
Since there exist extremal real numbers which are badly approximable
\cite{Ra}, this raises the question as to whether or not each
extremal real number is such.  At present, we simply know that an
extremal real number $\xi$ satisfies a measure of approximation by
rational numbers $p/q$ of the form
\begin{equation*}
 \Big| \xi - \frac{p}{q} \Big|
 \ge c_6 q^{-2} (1+\log |q|)^{-t},
\end{equation*}
with constants $c_6>0$ and $t\ge 0$ depending only on $\xi$ (Theorem
1.3 of \cite{Rb}).  In this paper, we establish a sufficient
condition for an extremal real number to have bounded partial
quotients and construct new examples of such numbers.

\section{Notation and statements of the main results}
 \label{sec:results}

A {\it Fibonacci sequence} in a monoid is a sequence $(w_i)_{i\ge
1}$ of elements of this monoid which satisfies the recurrence
relation $w_{i+2}=w_{i+1}w_i$ for each $i\ge 1$.  Here, we shall
work with two types of monoids.

One is the monoid of words $E^*$ on an alphabet $E$, with the
product given by concatenation of words.  A Fibonacci sequence
$(w_i)_{i\ge 1}$ in $E^*$ has the property that $w_i$ is a prefix
(left-factor) of $w_{i+1}$ for each $i\ge 2$ and so, it admits a
limit $w_\infty=\lim w_i$ in the completion of $E^*$ for pointwise
convergence. This limit is an infinite word unless $w_1$ and $w_2$
are empty.  For example, if $E=\{a,b\}$ consists of two distinct
elements $a$ and $b$, then the Fibonacci sequence of words starting
with $w_1=b$ and $w_2=a$ converges to the infinite word
$f_{a,b}=abaababa\dots$. In general, the limit of any Fibonacci
sequence of words $(w_i)_{i\ge 1}$ derives from this generic
infinite word $f_{a,b}$ by substituting into it the words $w_1$ and
$w_2$ for the letters $b$ and $a$ respectively.  It will come out
indirectly of our analysis that such a limit is an infinite
non-ultimately periodic word if and only if $w_1$ and $w_2$ do not
commute (see the remark after Theorem \ref{results:thmB} below).  A
direct proof of this fact has been recently provided by B.~Lucier
\cite{Lu}.

The other monoid is in fact a group.  It is constructed as follows.
Define the \emph{content} $c(A)$ of a non-zero matrix $A$ in $\MdZ$
to be the greatest positive common divisor of its coefficients and
say that such a matrix $A$ is \emph{primitive} if $c(A)=1$.  For
each non-zero $A\in\MdZ$, denote by $A^\red$ the unique primitive
integer matrix such that $A=c(A)A^\red$.  Then, the set $\cP$ of all
primitive matrices with non-zero determinant in $\MdZ$ is a group
for the operation $*$ given by $A*B = (AB)^\red$.  Its quotient
$\cP/\{\pm I\}$ is isomorphic to $\mathrm{PGL}_2(\bQ)$.

\begin{definition}
We say that a Fibonacci sequence $(W_i)_{i\ge 1}$ in $\cP$ is
\emph{admissible} if there exists a non-symmetric and
non-skew-symmetric matrix $N\in\cP$ such that, upon putting
$N_i=\tN$ for $i$ odd and $N_i=N$ for $i$ even, the product $W_iN_i$
is a symmetric matrix for each $i\ge 1$.
\end{definition}

This definition differs slightly from that in \S3 of \cite{Re}.
However, the same argument as in the proof of Proposition 3.1 of
\cite{Re} shows that most Fibonacci sequences in $\cP$ are
admissible in the sense that there exists a non-empty Zariski open
subset $\cU$ of $\GL_2(\bC)^2$ such that any pair $(W_1,W_2) \in
\cU\cap\cP$ generates an admissible Fibonacci sequence in $\cP$.

\medskip
We define also the {\it norm} $\|A\|$ of a matrix $A$ with real
coefficients to be the largest absolute value of its coefficients.
With this notation, we will prove in \S\ref{sec:thmA} the following
characterization of extremal real numbers which translates in the
present setting several results from \cite{Rb} and \cite{Rd}.

\begin{theorem}
 \label{results:thmA}
For each extremal real number $\xi$ there exists an unbounded
admissible Fibonacci sequence $(W_i)_{i\ge 1}$ in $\cP$ which
satisfies
\begin{equation}
 \label{thmA:eq1}
 \|W_{i+1}\| \ggll \|W_i\|^\gamma, \quad
 \|(\xi,-1)W_i\| \ggll \|W_i\|^{-1}
 \et
 |\det W_i| \ggll 1,
\end{equation}
with implied constants that are independent of $i$.  Such a sequence
is uniquely determined by $\xi$ up to its first terms, and up to
term-by-term multiplication by a Fibonacci sequence in $\{\pm 1\}$.
Conversely, any unbounded admissible Fibonacci sequence $(W_i)_{i\ge
1}$ in $\cP$ which satisfies
\begin{equation}
 \label{thmA:eq2}
 \|W_{i+2}\| \gg \|W_{i+1}\| \|W_i\|
 \et
 |\det W_i| \ll 1,
\end{equation}
also satisfies the conditions \eqref{thmA:eq1} for some extremal
real number $\xi$.
\end{theorem}

Thus, any unbounded admissible Fibonacci sequence in $\cP$
satisfying \eqref{thmA:eq2} is \emph{associated} to some extremal
real number $\xi$ in the sense that it satisfies \eqref{thmA:eq1}.
Note also that, since $\gamma^2=\gamma+1$, the first condition in
\eqref{thmA:eq1} is stronger than the first condition in
\eqref{thmA:eq2}.

It is shown in \S2 of \cite{Ra} and in \S6 of \cite {Rb} that, for
any choice of distinct positive integers $a$ and $b$, the real
number $\xi_{a,b}$ whose continued fraction expansion $\xi_{a,b} =
[0,a,b,a,a,\dots]$ is given by $0$ followed by the elements of
$f_{a,b}$ is an extremal real number.  More generally, we will prove
the following result (see \S\ref{sec:thmB}).

\begin{theorem}
 \label{results:thmB}
A real number $\xi$ is extremal with an associated Fibonacci
sequence in $\GL_2(\bZ)$ if and only if the sequence of its
partial quotients in its continued fraction expansion
coincides, up to its first terms, with the limit of a Fibonacci
sequence of words $(w_i)_{i\ge 1}$ in $(\bN \setminus
\{0\})^*$ starting with two non-commuting words $w_1$ and $w_2$.
\end{theorem}

Let $w_\infty = a_1 a_2 a_3 \dots$ be the limit of a Fibonacci
sequence of words $(w_i)_{i\ge 1}$ in $(\bN \setminus \{0\})^*$
starting with non-empty words $w_1$ and $w_2$. If $w_1$ and $w_2$
commute, then $w_\infty=\lim_{i\to\infty}(w_1)^i$ is a periodic word
and so $\xi=[0,a_1,a_2,\dots]$ is a quadratic real number.
Conversely, if $w_1$ and $w_2$ do not commute, the above theorem
shows that this real number $\xi$ is extremal. Since an extremal
real number is not quadratic, the infinite word $w_\infty$ cannot in
this case be ultimately periodic.

The next result provides a sufficient condition for an extremal
real number to be badly approximable.

\begin{theorem}
 \label{results:thmC}
Let $E=\{a,b\}$ be an alphabet of two letters, let $(w_k)_{k\ge 1}$
be the Fibonacci sequence in $E^*$ generated by $w_1=b$ and $w_2=a$,
and let $f_{a,b}=\lim_{k\to\infty}w_k$.  Let $\xi$ be an extremal
real number and let $(W_k)_{k\ge 1}$ be a Fibonacci sequence in $\cP$
which is associated to $\xi$.  Consider the morphism of
monoids $\Phi\colon E^* \to \cP$ mapping $w_k$ to $W_k$ for each
$k\ge 1$.  For each $i\ge 1$, denote by $u_i$ the prefix of $f_{a,b}$
with length $i$, and put $U_i=\Phi(u_i)$.  Then, we have
\begin{equation}
 \label{eq:Ui}
 \|(\xi,-1)U_i\| \ggll \frac{|\det U_i|}{\|U_i\|}
\end{equation}
with implied constants that do not depend on $i$.  Moreover, if the
sequence $(\det U_i)_{i\ge 1}$ is bounded, then $\xi$ is badly
approximable.
\end{theorem}

It would be interesting to know if, conversely, the sequence $(\det
U_i)_{i\ge 1}$ is bounded when $\xi$ is badly approximable.  The
proof of the above result is given in \S\ref{sec:thmC}.

Going back to the definitions, we note that, if $\xi$ is badly
approximable (resp.\ extremal) and if $a,b\in\bQ$ with $a\neq 0$,
then $a\xi+b$ and $1/\xi$ are as well badly approximable (resp.\
extremal). This implies that the set of badly approximable real
numbers is stable under the action of $\GL_2(\bQ)$ on
$\bR\setminus\bQ$ by linear fractional transformations.  Our last
main result, proved in \S\ref{sec:thmD}, is that there exist orbits
which do not contain any of the numbers produced by Theorem
\ref{results:thmB}.

\begin{theorem}
 \label{results:thmD}
There exist badly approximable extremal real numbers which are not
conjugate under the action of $\GL_2(\bQ)$ to any extremal real
number having an associated Fibonacci sequence in $\GL_2(\bZ)$.
\end{theorem}

\section{Proof of Theorem \ref{results:thmA}}
 \label{sec:thmA}

The following lemma gathers essentially all facts that we will need
from \cite{Rb} and \cite{Rd}.

\begin{lemma}
 \label{thmA:lemma1}
Let $\xi$ be an extremal real number.  Then, there exists an
unbounded sequence of symmetric matrices $(\uy_i)_{i\ge 1}$ in $\cP$
such that, for each $i\ge 1$, we have
\begin{equation}
 \label{thmA:lemma1:eq1}
 \|\uy_{i+1}\| \ggll \|\uy_i\|^\gamma, \quad
 \|(\xi,-1)\uy_i\| \ggll \|\uy_i\|^{-1}
 \et
 |\det \uy_i| \ggll 1,
\end{equation}
with implied constants that are independent of $i$. Such a sequence
$(\uy_i)_{i\ge 1}$ is uniquely determined by $\xi$ up to its first
terms and up to multiplication of each of its terms by $\pm 1$.
Moreover, for any such sequence, there exists a non-symmetric and
non-skew-symmetric matrix $M\in\cP$ such that
\begin{equation}
 \label{thmA:lemma1:eq2}
 \uy_{i+2}
 = \pm
   \begin{cases}
     \uy_{i+1}*M*\uy_i &\text{if $i$ is odd,}\\
     \uy_{i+1}*\tM*\uy_i &\text{if $i$ is even,}
   \end{cases}
\end{equation}
for any sufficiently large index $i$.  Conversely, if $(\uy_i)_{i\ge
1}$ is an unbounded sequence of symmetric matrices in $\cP$ which
satisfies a recurrence relation of the type \eqref{thmA:lemma1:eq2}
for some non-symmetric matrix $M\in\cP$, and if
\begin{equation}
 \label{thmA:lemma1:eq3}
 \|\uy_{i+2}\| \gg \|\uy_{i+1}\| \|\uy_i\|
 \et
 |\det \uy_i| \ll 1,
\end{equation}
then $(\uy_i)_{i\ge 1}$ also satisfies the estimates
\eqref{thmA:lemma1:eq1} for some extremal real number $\xi$.
\end{lemma}

\begin{proof}
The first assertion in this proposition comes from Theorem 5.1 of
\cite{Rb} upon noting that, for an arbitrary symmetric matrix $\uy =
\begin{pmatrix} y_0 &y_1\\ y_1 &y_2 \end{pmatrix}$, we have
\begin{equation*}
 \|(\xi,-1)\uy\|
 = \max\{ |y_0\xi-y_1|, |y_1\xi-y_2| \}
 \gg\ll \max\{ |y_0\xi-y_1|, |y_0\xi^2-y_2|\},
\end{equation*}
with implied constants depending only on $\xi$.  The second
assertion follows from Proposition 4.1 of \cite{Rd}, the third one
from Corollary 4.3 of \cite{Rd}, and the last one from Proposition
5.1 of \cite{Rd}.
\end{proof}

In the proof of Theorem \ref{results:thmA} below, we use repeatedly
the following observation (the proof of which is omitted).

\begin{lemma}
 \label{thmA:lemma2}
Let $(W_i)_{i\ge 1}$, $(\uy_i)_{i\ge 1}$ and $(N_i)_{i\ge 1}$ be
sequences in $\cP$, and let $\xi\in\bR$.  Assume that the sequence
$(N_i)_{i\ge 1}$ is bounded and that $\uy_i=W_i*N_i$ for each $i\ge
1$.  Then, we have
\begin{equation*}
 \|W_i\| \ggll \|\uy_i\|,
 \quad
 \|(\xi,-1)W_i\| \ggll \|(\xi,-1)\uy_i\|
 \et
 |\det(W_i)| \ggll |\det(\uy_i)|,
\end{equation*}
with implied constants that do not depend on $i$.
\end{lemma}

\medskip
\begin{proof}[\bf Proof of Theorem \ref{results:thmA}]
Let $\xi\in\bR$ be extremal.  Then, lemma \ref{thmA:lemma1} provides
an unbounded sequence of symmetric matrices $(\uy_i)_{i\ge 1}$ in
$\cP$ and a non-symmetric and non-skew-symmetric matrix $M\in\cP$
satisfying both the estimates \eqref{thmA:lemma1:eq1} and the
recurrence relation \eqref{thmA:lemma1:eq2} for each sufficiently
large $i$. Omitting if necessary a finite even number of initial
terms in the sequence $(\uy_i)_{i\ge 1}$, we may assume, without
loss of generality, that \eqref{thmA:lemma1:eq2} holds for each
$i\ge 1$. Then, for a suitable choice of signs, the formula
$W_i=\pm(\uy_i*M_i)$ with $M_i=\tM$ if $i$ is odd and $M_i=M$ if $i$
is even, defines an admissible Fibonacci sequence in $\cP$ (the
corresponding matrix $N$ is the inverse of $M$ in the group $\cP$).
Moreover, the estimates \eqref{thmA:lemma1:eq1} together with Lemma
\ref{thmA:lemma2} show that this sequence satisfies the conditions
\eqref{thmA:eq1} of Theorem \ref{results:thmA}.  This proves the
first assertion of the theorem.

Now, let $(W'_i)_{i\ge 1}$ be any unbounded admissible Fibonacci
sequence satisfying, like $(W_i)_{i\ge 1}$, the conditions
\eqref{thmA:eq1}, and let $N'\in\cP$ such that, upon putting
$N'_i=\tN'$ for $i$ odd and $N'_i=N'$ for $i$ even, the matrix
$\uy'_i=W'_i*N'_i$ is symmetric for each $i\ge 1$.  Then, lemma
\ref{thmA:lemma2} shows that $(\uy'_i)_{i\ge 1}$ satisfies, like
$(\uy_i)_{i\ge 1}$, the estimates \eqref{thmA:lemma1:eq1}.
Consequently, by Lemma \ref{thmA:lemma1}, there exist integers
$k,\ell\ge 0$ such that $\uy'_{i+k}=\pm\uy_{i+\ell}$ for each $i\ge
1$. Since we have $\uy_{i+2} = \pm(W_{i+1}*\uy_i)$ and $\uy'_{i+2} =
\pm(W'_{i+1}*\uy'_i)$ for $i\ge 1$, this implies that $W'_{i+k} =
\pm W_{i+\ell}$ for each $i\ge 2$.  Moreover, the signs $\pm$ in the
last formula must come from a Fibonacci sequence in $\{\pm 1\}$.
This proves the second assertion of the theorem.

Finally, let $(W_i)_{i\ge 1}$ be any unbounded admissible Fibonacci
sequence satisfying the conditions \eqref{thmA:eq2} in Theorem
\ref{results:thmA}, without reference to a given extremal real
number $\xi$, and let $N\in\cP$ such that, upon putting $N_i=\tN$
for $i$ odd and $N_i=N$ for $i$ even, the matrix $\uy_i=W_i*N_i$ is
symmetric for each $i\ge 1$.  Then, lemma \ref{thmA:lemma2} shows
that $(\uy_i)_{i\ge 1}$ satisfies the conditions
\eqref{thmA:lemma1:eq3} in Lemma \ref{thmA:lemma1}. Consequently it
satisfies the stronger conditions \eqref{thmA:lemma1:eq1} for some
extremal real number $\xi$, and thus, by Lemma \ref{thmA:lemma2},
satisfies the estimates \eqref{thmA:eq1} of Theorem
\ref{results:thmA} for the same $\xi$.
\end{proof}

\section{Proof of Theorem \ref{results:thmB}}
 \label{sec:thmB}

Serret's theorem asserts that two real numbers have continued
fraction expansions which coincide up to their first terms if and
only if these numbers belong to the same orbit under the action of
$\GL_2(\bZ)$ by linear fractional transformations (Theorem 6B of
\cite{Sc}).  Our proof for Theorem \ref{results:thmB} is inspired
from the proof of this result given by Cassels in \S3, Chap.\ I of
\cite{Ca}. We break it into two propositions.  To establish the
first one, we need the following auxiliary result which provides the
link with continued fractions.

\begin{lemma}
 \label{B:lemma1}
Let $A=\begin{pmatrix} a&b\\ c&d\end{pmatrix}\in \GL_2(\bZ)$ with
$d\ge 1$.  Then, there is one and only one choice of integers $s\ge
1$ and $a_0,a_1,\dots,a_s$ with $a_1,\dots,a_{s-1}\ge 1$ such that
\begin{equation}
 \label{B:eqA}
 A = \begin{pmatrix} a_0&1\\ 1&0\end{pmatrix}
     \begin{pmatrix} a_1&1\\ 1&0\end{pmatrix}
     \cdots
     \begin{pmatrix} a_s&1\\ 1&0\end{pmatrix}.
\end{equation}
These integers are also characterized by the properties
\begin{equation}
 \label{B:cfrac}
 \frac{b}{d}=[a_0,\dots,a_{s-1}], \quad
 \frac{c}{d}=[a_s,\dots,a_1], \quad
 \det(A)=(-1)^{s+1}.
\end{equation}
\end{lemma}

\begin{proof}
Induction on $s$ shows that, if $A$ can be written in the form
\eqref{B:eqA} for a choice of integers $s \ge 1$ and $a_0, a_1,
\dots, a_s$ with $a_1, \dots, a_{s-1} \ge 1$, then we have $b/d =
[a_0, \dots, a_{s-1}]$.  Taking the transpose of both sides of
\eqref{B:eqA}, this observation also provides $c/d = [a_s, \dots,
a_1]$. Moreover the last equality in \eqref{B:cfrac} follows from
the multiplicativity of the determinant.  Since each rational number
has exactly two continued fraction expansions with lengths differing
by one, this proves the uniqueness of the factorization
\eqref{B:eqA}, when it exists.

Now, without making assumptions on $A$, define an integer $s\ge 1$
and a sequence of integers $a_1, \dots, a_s$ with $a_1, \dots,
a_{s-1} \ge 1$ by the conditions $c/d = [a_s, \dots, a_1]$ and
$\det(A) = (-1)^{s+1}$.  Define also $a_0$ to be the integer for
which the distance between $b/d$ and $[a_0,\dots,a_{s-1}]$ is at
most $1/2$. Then, by the above observations, the right hand side of
\eqref{B:eqA} is a matrix $A' =
\begin{pmatrix} a'&b'\\ c'&d'\end{pmatrix}$ with the same
determinant as $A$, satisfying $d'\ge 1$, $c'/d'=c/d$ and
$|b'/d'-b/d|\le 1/2$.  Since $(c,d)$ and $(c',d')$ are rows of
matrices in $\GL_2(\bZ)$, they are primitive points of $\bZ^2$ and
the relation $c'/d'=c/d$ implies $(c',d')=\pm (c,d)$.  Since $d$ and
$d'$ are positive, we deduce that $A$ and $A'$ have the same second
row $(c',d')=(c,d)$.  Since these matrices also have the same
determinant, this forces $(a',b')=(a,b)+k(c,d)$ for some integer
$k$. Then, we find $|b'/d'-b/d| = |k|$, thus $k=0$ and therefore
$A'=A$.
\end{proof}

\begin{corollary}
 \label{B:cor1}
Let $\cS_1$ denote the set of matrices $\begin{pmatrix} a&b\\
c&d\end{pmatrix}\in \MdR$ with $a\ge \max\{b,c\}$ and $\min\{b,c\}
\ge d\ge 0$, and define $\cS=\cS_1\cap \GL_2(\bZ)$.  Then, $\cS_1$
and $\cS$ are closed under multiplication and transposition.
Moreover, the map from $\bN\setminus\{0\}$ to $\cS$ sending $a$ to
$\begin{pmatrix}a&1\\1&0\end{pmatrix}$ for each
$a\in\bN\setminus\{0\}$ extends to an isomorphism of monoids
$\sigma\colon(\bN\setminus\{0\})^* \to \cS\cup\{I\}$.
\end{corollary}

\begin{proof}
The only delicate point here is the surjectivity of the map
$\sigma$.  Clearly, any matrix $A=\begin{pmatrix}a&b\\c&d
\end{pmatrix} \in\cS$ is in the image of $\sigma$ if $d=0$
because then we have $a\ge b=c=1$.  If $d\ge 1$, we note that
the integers $a_0$ and $a_s$ given by Lemma \ref{B:lemma1}
are respectively the integral parts of $b/d$ and $c/d$ except
in the case where $d=1$ and $\det A=-1$.  In the latter case, we
have $s=2$, $a_0=b-1$, $a_1=1$, $a_2=c-1$, and also $a=bc-1$.
Then, the condition $a\ge \max\{b,c\}\ge 1$ implies $b,c\ge
2$ and so $a_0,a_s\ge 1$.  Otherwise, the condition
$\min\{b,c\} \ge d\ge 1$ ensures that $a_0,a_s\ge 1$.  So, in both
cases, the integers $a_0,\dots,a_s$ are positive, and $A$ is the
image of $(a_0,\dots,a_s)$ under $\sigma$.
\end{proof}

The following proposition presents a first step towards the proof of
Theorem \ref{results:thmB}.

\begin{proposition}
 \label{B:prop1}
The set of extremal real numbers with an associated (admissible)
Fibonacci sequence in $\GL_2(\bZ)$ is stable under the action of
$\GL_2(\bZ)$ by linear fractional transformations.  Any orbit
contains an extremal real number with an associated Fibonacci
sequence in $\cS$.
\end{proposition}

\begin{proof}
Let $\xi$ be an extremal real number with an associated Fibonacci
sequence $(W_i)_{i\ge 1}$ in $\GL_2(\bZ)$.  This sequence being
admissible, there exists $N\in\cP$ with $N\neq \pm\tN$ such that,
upon putting $N_i=N$ if $i$ is even and $N_i=\tN$ if $i$ is odd, the
product $\uy_i = W_iN_i$ is symmetric for each $i\ge 1$.

For each $U\in\GL_2(\bZ)$, the sequence $(W'_i)_{i\ge 1} =
(U^{-1}W_iU)_{i\ge 1}$ is a Fibonacci sequence in $\GL_2(\bZ)$. It
is admissible with corresponding matrix $N'=U^{-1} N \tU^{-1}$, and
it satisfies the conditions \eqref{thmA:eq1} of Theorem
\ref{results:thmA} with $W_i$ replaced by $W'_i$ and $\xi$ replaced
by the real number $\eta$ such that $(\eta,-1)$ is proportional to
$(\xi,-1)U$.  By varying $U$, we get in this way all real numbers
$\eta$ which are conjugate to $\xi$ under $\GL_2(\bZ)$.  So, these
numbers are extremal with an associated Fibonacci sequence in
$\GL_2(\bZ)$. This proves the first assertion of the lemma.

For the second assertion, let $1/\xi=[a_0,a_1,a_2,\dots]$ be the
continued fraction expansion of $1/\xi$.  Put $d=\det(N)$ and
$M=dN^{-1}$, so that we have $M\in\cP$ and $W_i=d^{-1}\uy_iM_i$
where $M_i=M$ if $i$ is even and $M_i=\tM$ if $i$ is odd.  Since $M$
is not skew-symmetric and since $[\bQ(\xi):\bQ]>2$, the product
\begin{equation*}
 \theta
  = \begin{pmatrix} 1/\xi&1 \end{pmatrix}
    M
    \begin{pmatrix} 1/\xi\\1 \end{pmatrix}
\end{equation*}
is non-zero.  Replacing $N$ by $-N$ if necessary, so that $M$ is
replaced by $-M$, we may assume without loss of generality that this
number $\theta$ is positive.  For each $k\ge 1$, define
\begin{equation*}
 U_k =
     \begin{pmatrix} a_0&1\\ 1&0\end{pmatrix}
     \begin{pmatrix} a_1&1\\ 1&0\end{pmatrix}
     \cdots
     \begin{pmatrix} a_k&1\\ 1&0\end{pmatrix}.
\end{equation*}
Then, the standard recurrence relations in the theory of continued
fractions show that we have $U_k =
\begin{pmatrix} p_k &p_{k-1}\\ q_k &q_{k-1}\end{pmatrix}$
where $p_k/q_k=[a_0,\dots,a_k]$ denotes the $k$-th convergent of
$1/\xi$ written in reduced form.  Since $|q_k(1/\xi)-p_k|\le
1/q_{k+1}$ for each $k\ge 0$, this gives
\begin{equation*}
 U_k = \begin{pmatrix} 1/\xi\\ 1\end{pmatrix}
       \begin{pmatrix} q_k &q_{k-1}\end{pmatrix}
       + \cO(1/q_{k})
\end{equation*}
and thus
\begin{equation*}
 \tU_k M U_k
  = \theta
    \begin{pmatrix}
    q_k^2 &q_{k-1}q_k\\ q_{k-1}q_k &q_{k-1}^2
    \end{pmatrix}
    + \cO(1).
\end{equation*}
The latter matrix belongs to $\cS_1$ if $k$ sufficiently large,
because we have $q_k>q_{k-1}$ for each $k\ge 2$, and $q_{k-1}$ tends
to infinity with $k$.  Fix such a value of $k$.  Since $\cS_1$ is
closed under transposition, we get $\tU_kM_iU_k\in\cS_1$ for each
$i\ge 1$.  We claim that $\epsilon_i U_k^{-1}W_iU_k$ also belongs to
$\cS_1$ for an appropriate choice of $\epsilon_i\in\{-1,1\}$ and
each sufficiently large $i$. To prove this, we note that the product
$U_k^{-1}\begin{pmatrix}1/\xi\\1\end{pmatrix}$ is proportional to
$\begin{pmatrix}r\\1\end{pmatrix}$ where $r=[a_{k+1},a_{k+2},\dots]$
is a real number with $r>1$.  Since, for each $i\ge 1$, we have
\begin{equation*}
 \uy_i
  = y_{i,2}
    \begin{pmatrix} \xi^{-2} &\xi^{-1}\\ \xi^{-1} &1\end{pmatrix}
    + \cO(\|W_i\|^{-1})
\end{equation*}
with $y_{i,2}\in\bZ$, we find
\begin{equation*}
 U_k^{-1}\uy_i\tU_k^{-1} = c_i \begin{pmatrix} r^2 &r\\ r &1\end{pmatrix}
         + \cO(\|W_i\|^{-1}),
\end{equation*}
for some $c_i\in\bR$.  Thus, if $i$ is sufficiently large, say $i\ge
i_0$, the matrix $\pm U_k^{-1}\uy_i\tU_k^{-1}$ belongs to $\cS_1$
for an appropriate choice of sign $\pm$.  Multiplying this matrix on
the right by $\tU_kM_iU_k$ which also belongs to $\cS_1$, we deduce
that $\pm U_k^{-1} \uy_i M_i U_k \in \cS_1$ for the same choice of
sign and thus that $\epsilon_i U_k^{-1}W_iU_k\in\cS_1$ for some
$\epsilon_i\in\{-1,1\}$.  Since $U_k\in\GL_2(\bZ)$ and since
$W_i\in\GL_2(\bZ)$ for each $i\ge 1$, we conclude that $(\epsilon_i
U_k^{-1}W_iU_k)_{i\ge i_0}$ is an admissible Fibonacci sequence in
$\cS$.  By the first part of the proof, it is associated to an
extremal real number $\eta$ in the same $\GL_2(\bZ)$-orbit as $\xi$.
\end{proof}

We also need the following technical
result.

\begin{lemma}
 \label{B:lemma2}
Let $(W_i)_{i\ge 1}$ be a Fibonacci sequence in $\cP$.  If $W_1$ and
$W_2$ do not have a common eigenvector in $\bQ^2$ and satisfy
$W_1W_2\neq \pm W_2W_1$, then $(W_i)_{i\ge 1}$ is an admissible
Fibonacci sequence.
\end{lemma}

\begin{proof}
We first note that there exists a non-zero primitive matrix
$N\in\MdZ$ such that $W_1\tN$, $W_2N$ and $W_3\tN$ are symmetric
because these three conditions translate into a system of three
homogeneous linear equations in the four unknown coefficients of
$N$.  Fix such a choice of $N$ and define accordingly $N_i=\tN$ for
$i$ odd and $N_i=N$ for $i$ even.  Then, the product $W_iN_i$ is
symmetric for $i=1,2,3$ and using the relation of proportionality
\begin{equation*}
 W_{i+3}N_{i+3}
 \propto
 (W_{i+1}N_{i+1}) N_{i+1}^{-1} (W_iN_i) N_i^{-1} (W_{i+1}N_{i+1}),
\end{equation*}
we deduce by induction on $i$ that $W_iN_i$ is symmetric for each
$i\ge 1$.

If $\det N=0$, then we can write $N=A\tB$ with non-zero column
vectors $A$ and $B$ in $\bQ^2$.   Since $W_1\tN=(W_1B)\tA$ is
symmetric, we deduce that $W_1B\propto A$. Similarly, since
$W_2N=(W_2A)\tB$ and $W_3\tN=(W_3B)\tA$ are symmetric, we find that
$W_2A\propto B$ and $W_3B\propto A$.  Using the first two relations
of proportionality, we also get $W_3B\propto W_2(W_1B) \propto W_2A
\propto B$.  As $W_3B\neq 0$, this shows that $A\propto B$, and thus
that $B$ is a common eigenvector of $W_1$ and $W_2$, against the
hypothesis.  Thus we have $N\in\cP$.  We also note that
\begin{equation*}
 W_2W_1\tN
  = {^t}(W_2W_1\tN)
  = {^t}(W_1\tN)\tW_2
  = (W_1\tN)\tW_2
  = W_1 {^t}(W_2N)
  = W_1W_2 N.
\end{equation*}
Since $W_1W_2\neq \pm W_2W_1$, this implies that $\tN\neq \pm N$.
Thus, the sequence $(W_i)_{i\ge 1}$ is admissible.
\end{proof}

The hypotheses of Lemma \ref{B:lemma2} are satisfied for example
when the matrices $W_1$, $W_2$, $W_1W_2$ and $W_2W_1$ are linearly
independent over $\bQ$.  The corollary below provides another
instance where this lemma applies.

\begin{corollary}
 \label{B:cor2}
Any Fibonacci sequence $(W_i)_{i\ge 1}$ in $\cS$ generated by two
non-commuting matrices $W_1, W_2\in \cS$ is admissible.
\end{corollary}

\begin{proof}
Since $\cS\subset\GL_2(\bZ)$, the eigenvalues of a matrix $W\in\cS$
are algebraic units.  So, if one of them is rational, both of them
belong to $\{-1,1\}$.  Since the only matrices of $\cS$ with trace
at most $2$ are $\begin{pmatrix}1&1\\1&0\end{pmatrix}$ and
$\begin{pmatrix}2&1\\1&0\end{pmatrix}$ which have no rational
eigenvalue, we deduce that no matrix of $\cS$ has a rational
eigenvalue.  In particular, any $W_1, W_2\in \cS$ do not share a
common eigenvector in $\bQ^2$. Since such matrices have non-negative
coefficients and non-zero product, they also satisfy $W_1W_2\neq
-W_2W_1$.  Thus, if they do not commute, lemma \ref{B:lemma2} shows
that they generate an admissible Fibonacci sequence.
\end{proof}

Serret's theorem combined with Proposition \ref{B:prop1} reduces the
proof of Theorem \ref{results:thmB} to the following statement.

\begin{proposition}
 \label{B:prop2}
A real number $\xi$ is extremal with an associated Fibonacci
sequence in $\cS$ if and only if its continued fraction expansion
is of the form $[0,a_1,a_2,\dots]$ where $(a_1,a_2,\dots)$ is the
limit of a Fibonacci sequence of words $(w_i)_{i\ge 1}$ in $(\bN \setminus
\{0\})^*$ starting with two non-commuting words $w_1$ and $w_2$.
\end{proposition}

\begin{proof}
Let $\xi = [0,a_1,a_2,\dots]$ where $(a_1,a_2,\dots)$ is the limit
of a sequence of words $(w_i)_{i\ge 1}$ in $(\bN \setminus \{0\})^*$
starting with two non-commuting words $w_1$ and $w_2$. Denote by
$(W_i)_{i\ge 1}$ the image of the sequence $(w_i)_{i\ge 1}$ under
the isomorphism of monoids $\sigma\colon (\bN \setminus \{0\})^* \to
\cS\cup\{I\}$ defined in Corollary \ref{B:cor1}.  Since $w_1$ and
$w_2$ do not commute, the same is true of $W_1$ and $W_2$ and so, by
Corollary \ref{B:cor2}, $(W_i)_{i\ge 1}$ is an admissible Fibonacci
sequence in $\cS$. We also note that, for each pair of matrices
$A,B\in\cS$, we have $\|AB\| > \|A\| \|B\|$.  Then, the relation
$W_{i+2}=W_{i+1}W_i$ implies $\|W_{i+2}\| > \|W_{i+1}\| \|W_i\|$ for
each $i\ge 1$. In particular, the sequence $(W_i)_{i\ge 1}$ is
unbounded.  As $|\det W_i|=1$ for each $i$, it also satisfies the
conditions \eqref{thmA:eq2} of Theorem \ref{results:thmA}. Thus, the
sequence $(W_i)_{i\ge 1}$ is associated to some extremal real number
$\eta$.   On the other hand, the theory of continued fractions shows
that $\|(\xi,-1)W_i\| \ll \|W_i\|^{-1}$ since the ratios of the
elements in the columns of $W_i$ are successive convergents of
$1/\xi$.  Thus, $\xi=\eta$ is extremal.

Conversely, let $\xi$ be an extremal real number with an associated
Fibonacci sequence $(W_i)_{i\ge 1}$ in $\cS$.  The inverse image
of this sequence under $\sigma$ is a Fibonacci sequence
$(w_i)_{i\ge 1}$ in $(\bN \setminus \{0\})^*$ and, as above, we deduce
that $\xi=[0,a_1,a_2,\dots]$ where $(a_1,a_2,\dots) =
\lim_{i\to\infty} w_i$.  Since $\xi$ is neither rational nor quadratic,
this sequence is infinite and ultimately not periodic.  In particular,
$w_1$ and $w_2$ are not both powers of the same word, and so they do
not commute (Proposition 1.3.2 of Chapter 1 of \cite{Lo}).
\end{proof}

\begin{remark}
Let $w_\infty=(1,2,3,1,2,1,2,3,\dots)$ be the limit of the Fibonacci
sequence $(w_i)_{i\ge 1}$ generated by $w_1=(3)$ and $w_2=(1,2)$.
Since $w_1w_2\neq w_2w_1$, the corresponding real number $\xi =
[0,1,2,3,1,2,1,2,3,\dots]$ is extremal.  However, contrary to the
generic Fibonacci word $f_{a,b}$ which contains palindromes of
arbitrary length as prefixes, the infinite word $w_\infty$ contains
no factor of length greater than $3$ which is a palindrome.
\end{remark}

\section{Proof of Theorem \ref{results:thmC}}
 \label{sec:thmC}

Throughout this section, the notation is the same as in Theorem
\ref{results:thmC}. Namely, we fix an alphabet $E=\{a,b\}$ of two
letters and denote by $(w_k)_{k\ge 1}$ the Fibonacci sequence in
$E^*$ generated by $w_1=b$ and $w_2=a$, with limit $f_{a,b}$.  We
also fix an extremal real number $\xi$ with an associated Fibonacci
sequence $(W_k)_{k\ge 1}$ in $\cP$, and denote by $\Phi\colon E^*
\to \cP$ the morphism of monoids mapping $w_k$ to $W_k$ for each
$k\ge 1$.  We start with the following observation.

\begin{lemma}
 \label{lemma:www}
Let $k$ and $\ell$ be integers with $k\ge \ell\ge 2$, and let
$w_k=uv$ be a factorization of $w_k$ in $E^*$.  Then, there exist a
prefix $u_0$ of $w_\ell$ and strictly decreasing sequences of
integers $i_1>i_2>\cdots>i_s$ and $j_1>j_2>\cdots>j_t$ bounded below
by $\ell$ such that
\begin{equation*}
 u = w_{i_1} w_{i_2} \cdots w_{i_s}u_0
 \et
 u_0v = w_{j_t} \cdots w_{j_2} w_{j_1}.
\end{equation*}
If $u$ is not a prefix of $w_\ell$, we can ask that $i_1\le k-1$ and
$j_1\le k-2$.
\end{lemma}

\begin{proof}
If $u$ is a prefix of $w_\ell$, we take $u_0=u$ so that $u_0v=w_k$.
Otherwise, we have $k>\ell$, thus $k\ge 3$ and the factorization
$w_k=w_{k-1}w_{k-2}$ implies that either there is a word $u'$ such
that $u=w_{k-1}u'$ and $u'v=w_{k-2}$, or we have $k\ge \ell+2$ and
there is a word $v'$ such that $v=v'w_{k-2}$ and $uv'=w_{k-1}$. The
result then follows by induction on $k$.
\end{proof}

Since the sequence $(W_i)_{i\ge 1}$ is admissible, there exists a
non-symmetric and non-skew-symmetric matrix $N$ such that, upon
putting $N_i=N$ if $i$ is even and $N_i=\tN$ if $i$ is odd, the
product $\uy_i=W_iN_i$ is symmetric for each $i\ge 1$.  This matrix
$\uy_i$ may not be primitive but, for the next result, it is
convenient not to normalize it.

\begin{lemma}
 \label{lemma:ymym}
Define $L = \max\{1,|\xi|\}^{-1} (1,\xi)$ and $\theta = L N^{-1}
(\tL)$. Then, there exist an index $\ell\ge 1$ and a constant $c\ge
1$ such that, for any sequence of integers $(i_1,\dots,i_s)$ with
entries bounded below by $\ell$ and repeated at most twice, we have
\begin{equation*}
 \frac{1}{c}
 \le
 \frac{
 \|W_{i_1}W_{i_2}\cdots W_{i_s}\|}
 { |\theta|^s \|\uy_{i_1}\|\|\uy_{i_2}\|\cdots \|\uy_{i_s}\|}
 \le
 c.
\end{equation*}
\end{lemma}

Note that we have $\theta\neq 0$ since $\xi$ is transcendental and
$N$ is not skew-symmetric.

\begin{proof}
Write $\uy_i = \begin{pmatrix} y_{i,0} &y_{i,1}\\ y_{i,1} &y_{i,2}
\end{pmatrix}$ for each $i\ge 1$.  As $\|(\xi,-1)\uy_i\|\ll
\|(\xi,-1)W_i\| \ll \|W_i\|^{-1}$, we have
\begin{equation*}
 \uy_i
  = y_{i,0}
    \begin{pmatrix} 1\\ \xi \end{pmatrix}
    \begin{pmatrix} 1  &\xi \end{pmatrix}
    +
    \cO\big(\|W_i\|^{-1}\big),
\end{equation*}
and so $\|\uy_i\| = |y_{i,0}| \max\{1,|\xi|\}^2 + \cO\big(
\|W_i\|^{-1} \big)$.  In particular, this shows that $y_{i,0}\neq 0$
for each sufficiently large $i$, say for $i\ge \ell$.  Then, for
those values of $i$, we find
\begin{equation}
 \label{eq:est:YM}
 \frac{W_i}{\theta \|\uy_i\|}
 =
 A_i + R_i
\end{equation}
where $R_i = \cO\big(\|W_i\|^{-2}\big)$ and where $A_i = \pm
\theta^{-1} (\tL) L N_i^{-1}$ belongs to the set
\begin{equation*}
 \cA
  = \Big\{ \pm I,\ \pm\frac{1}{\theta} \tL L N^{-1},\
           \pm\frac{1}{\theta} \tL L \tN^{-1} \Big\}.
\end{equation*}
Since $\theta = L N^{-1} (\tL) = L (\tN^{-1}) (\tL)$, the set $\cA$
is stable under multiplication.

Now, let $(i_1,\dots,i_s)$ be any sequence of integers bounded below
by $\ell$, with no entry repeated more than twice.  Using
\eqref{eq:est:YM}, we find
\begin{equation*}
 \frac{W_{i_1}\cdots W_{i_s}}
      {\theta^s \|\uy_{i_1}\| \cdots \|\uy_{i_s}\|}
 = A + R
\end{equation*}
where $A = A_{i_1}\cdots A_{i_s}$ belongs to $\cA$ and where $R$ is
a sum, indexed by all non-empty subsequences $(j_1,\dots,j_t)$ of
$(i_1,\dots,i_s)$, of products of the form $B_1R_{j_1}\cdots
B_tR_{j_t}B_{t+1}$ with $B_1,\dots,B_{t+1}\in\cA$.  Thus, for an
appropriate constant $\kappa>0$, we have
\begin{equation*}
 \|R\|
 \le \big(1+ \kappa\|W_{i_1}\|^{-2}\big) \cdots
               \big(1+ \kappa\|W_{i_s}\|^{-2}\big) - 1
 \le \exp\Big(2\kappa\sum_{i=\ell}^\infty \|W_i\|^{-2}\Big) - 1.
\end{equation*}
If $\ell$ is sufficiently large, this gives $\|R\| \le \|A\|/2$, and
so $\|A+R\| \ggll \|A\|\ggll 1$, as requested.
\end{proof}

\begin{lemma}
 \label{lemma:UV}
Let $k$ be a positive integer and let $w_k=uv$ be a factorization of
$w_k$ in $E^*$.  Put $U=\Phi(u)$ and $V=\Phi(v)$.  Then, we have
$\|UV\|\ggll \|U\|\|V\|$ with implied constants that are independent
of $k$, $u$ and $v$.
\end{lemma}

\begin{proof}
Let $\ell$ be as in Lemma \ref{lemma:ymym}. Without loss of
generality, we may assume that $k\ge \ell$.  Then, according to
Lemma \ref{lemma:www}, we can write $u = w_{i_1}\cdots w_{i_s}u_0$
and $u_0v = w_{j_t} \cdots w_{j_1}$ where $u_0$ is a prefix of
$w_\ell$ and where $(i_1,\dots,i_s)$ and $(j_1,\dots,j_t)$ are
strictly decreasing sequences of integers bounded below by $\ell$.
Put
\begin{equation*}
 P = W_{i_1}\cdots W_{i_s},
 \quad
 Q = W_{j_t}\cdots W_{j_1}
 \et
 U_0 = \Phi(u_0).
\end{equation*}
Then, we have $U=aPU_0$ and $U_0V=bQ$ with non-zero rational numbers
$a$ and $b$.  Since $U_0$ belongs to a finite set of matrices in
$\cP$, we deduce that
\begin{equation*}
 \|U\| \ggll |a| \|P\|
 \et
 \|V\| \ggll |b| \|Q\|.
\end{equation*}
Moreover, since the sequence $(i_1,\dots,i_s,j_t,\dots,j_1)$ has its
entries repeated at most twice and bounded below by $\ell$, Lemma
\ref{lemma:ymym} gives
\begin{equation*}
 \|PQ\|
 \ggll
 |\theta|^{s+t}
 \|\uy_{i_1}\| \cdots \|\uy_{i_s}\| \|\uy_{j_t}\| \cdots \|\uy_{j_1}\|
 \ggll
 \|P\|\|Q\|.
\end{equation*}
The conclusion follows because $UV = ab PQ$.
\end{proof}

\medskip
\begin{proof}[\bf Proof of Theorem \ref{results:thmC}]
We first note that, for any $(x,y)\in\bR^2$ and $U\in\GL_2(\bR)$, we
have
\begin{equation}
 \label{eq:est:U}
 \|(x,y)U\|
  \ge \frac{\|(x,y)UU^{-1}\|}{2\|U^{-1}\|}
   =  \frac{\|(x,y)\|\,|\det U|}{2\|U\|}.
\end{equation}
Applying this to the point $(x,y) = (\xi,-1)$ and the matrix $U_i =
\Phi(u_i)$ where $u_i$ denotes the prefix of $f_{a,b}$ of length
$i$, we get
\begin{equation*}
 \|(\xi,-1)U_i\| \ge \frac{|\det U_i|}{2\|U_i\|}
\end{equation*}
for each $i\ge 1$.  To prove an upper bound of the same type for
$\|(\xi,-1)U_i\|$, we denote by $k=k(i)$ the smallest positive
integer such that $u_i$ is a prefix of $w_k$, and write $w_k=u_iv_i$
with $v_i\in E^*$.  Putting $V_i=\Phi(v_i)$, we then have
\begin{equation}
 \label{eq:est:UV}
 U_iV_i = m_i W_k
\end{equation}
for some integer $m_i\ge 1$.  Applying \eqref{eq:est:U} to the point
$(x,y) = (\xi,-1)U_i$ and the matrix $U=V_i$, we find
\begin{equation*}
 \|(\xi,-1)W_k\|
  = \frac{1}{m_i} \|(\xi,-1)U_iV_i\|
  \ge \frac{\|(\xi,-1)U_i\|\, |\det V_i|}{2m_i\|V_i\|}.
\end{equation*}
Since $\|(\xi,-1)W_k\| \ll \|W_k\|^{-1}$, this gives
\begin{equation}
 \label{eq:est:ub1}
 \|(\xi,-1)U_i\|\
  \ll \frac{m_i \|V_i\|}{\|W_k\|\,|\det V_i|}.
\end{equation}
Applying Lemma \ref{lemma:UV} to the factorization
\eqref{eq:est:UV} on one hand, and taking determinants of both sides
of \eqref{eq:est:UV} on the other hand, we also find
\begin{equation*}
 \|V_i\| \ggll \frac{m_i \|W_k\|}{\|U_i\|}
 \et
 |\det V_i|
  = \frac{m_i^2 |\det W_k|}{|\det U_i|}
  \ge \frac{m_i^2}{|\det U_i|}.
\end{equation*}
These estimates combined with \eqref{eq:est:ub1} lead to
\begin{equation}
 \label{eq:est:ub2}
 \|(\xi,-1)U_i\|\ \ll \frac{|\det U_i|}{\|U_i\|},
\end{equation}
which completes the proof of \eqref{eq:Ui} in Theorem
\ref{results:thmC}.

Now, assume that the integers $\det U_i$ are bounded independently
of $i$ and, for each $i\ge 1$, choose a column $\begin{pmatrix} q_i\\
p_i \end{pmatrix}$ of $U_i$ with the largest norm.  Then,
\eqref{eq:est:ub2} leads to
\begin{equation*}
 |q_i\xi-p_i| \ll \|U_i\|^{-1}.
\end{equation*}
Since $U_{i+1}$ is either equal to $U_i*W_1$ or to $U_i*W_2$, we
also have $\|U_{i+1}\| \ll \|U_i\|$ and thus $|q_{i+1}| \ll \|U_i\|
\ll |q_i|$. Combining these estimates and noting that
$\gcd(p_i,q_i)$ is a divisor of $\det U_i$, we deduce the existence
of a constant $c\ge 1$ such that
\begin{equation}
 \label{eq:est:ub3}
 |q_{i+1}(q_i\xi-p_i)| \le c
 \et
 |q_{i+1}|\le c \Big| \frac{q_i}{\gcd(p_i,q_i)} \Big|,
\end{equation}
for each $i\ge 1$.  Moreover, we have $\limsup_{i\to \infty} |q_i| =
\infty$ since $(U_i)_{i\ge 1}$ contains the unbounded sequence
$(W_k)_{k\ge 1}$ as a subsequence. These facts imply that $\xi$ is
badly approximable. Indeed, if $p/q$ is an arbitrary rational
number, then, at the expense of replacing $c$ by a larger constant
if necessary, we may assume that there exists an index $i\ge 1$ such
that $ 0 < |q_i| \le 2c|q| \le |q_{i+1}|$. Using \eqref{eq:est:ub3},
this gives $2|q| \le |q_i/\gcd(p_i,q_i)|$, thus $p/q\neq p_i/q_i$
and so we find
\begin{equation*}
 \Big| \xi-\frac{p}{q} \Big|
  \ge \Big| \frac{p_i}{q_i}-\frac{p}{q} \Big|
     - \Big| \frac{p_i}{q_i}-\xi \Big|
  \ge \frac{1}{|qq_i|} -\frac{c}{|q_iq_{i+1}|}
  \ge \frac{1}{2|qq_i|}
  \ge \frac{1}{4cq^2}.
\end{equation*}
\end{proof}

\section{Proof of Theorem \ref{results:thmD}}
 \label{sec:thmD}

Again, let $E$ be a set of two elements $a$ and $b$, and let
$(w_i)_{i\ge 1}$ be the Fibonacci sequence in $E^*$ determined by
the conditions $w_1=b$ and $w_2=a$, with limit $f_{a,b}$.  The
following lemma is our main-tool for constructing more extremal real
numbers.

\begin{lemma}
 \label{D:lemma1}
Let $m$ be a non-zero integer and let $W\in \MdZ$ with $W^2 \equiv
0$ mod $m$.  Assume that there exist primitive matrices $W_1, W_2
\in \MdZ$ of determinant $m$ with $W_1\equiv W_2 \equiv W$ mod $m$,
and consider the morphism of monoids $\Phi\colon E^* \to \cP$
mapping $a$ to $W_2$ and $b$ to $W_1$.  Then, for each word $u\in
E^*$, the determinant of $\Phi(u)$ is $1$ if $u$ has even length and
it is $m$ if $u$ has odd length.
\end{lemma}

\begin{proof}
We proceed by recurrence on the length $\ell$ of $u$.  If $\ell \le
1$, the result is clear (for the empty word $1$, the matrix
$\Phi(1)$ is the identity).   If $\ell = 2$, we have
$\Phi(u)=(W_iW_j)^\red$ for some choice of indices $i,j \in
\{1,2\}$.  Then, since $W_iW_j \equiv W^2 \equiv 0$ mod $m$ and
since $\det(W_iW_j)=m^2$, the matrix $W_iW_j$ has content $|m|$, and
so $\Phi(u) = |m|^{-1} W_iW_j$ has determinant $1$.  Now, assume
that $\ell>2$ and that the result is true for words of smaller
length.  Write $u=u'u''$ where $u'$ has even length and $u''$ has
length $1$ or $2$.  By induction hypothesis, $\Phi(u')$ has
determinant $1$ while $\Phi(u'')$ is primitive with $\det \Phi(u'')
= 1$ if $\ell$ is even and $\det \Phi(u'') = m$ if $\ell$ is odd.
Then the product $\Phi(u')\Phi(u'')$ is primitive and so $\Phi(u) =
\Phi(u') \Phi(u'')$ has the same determinant as $\Phi(u'')$.
\end{proof}

We also need the following technical result.

\begin{lemma}
 \label{D:lemma2}
Let $r$ be a real number with $0<r\le 1$ and let $\cS_r$ denote the
set of matrices $A\in\MdR$ with positive coefficients whose elements
of the first row are bounded below by $r$ times those of the second
row, and whose elements of the first column are bounded below by $r$
times those of the second column. Then, $\cS_r$ is closed under
multiplication and, for each $A,A'\in\cS_r$, we have $\|AA'\| >
r\|A\|\|A'\|$.
\end{lemma}

\begin{proof}
The set $\cS_r$ consists of all $2\times 2$ matrices $A$ with
positive coefficients such that the products $(1, -r)A$ and
$(1,-r)\tA$ have non-negative coefficients.  The fact that this set
is closed under multiplication then follows from the associativity
of the matrix product.  To prove the second assertion, take
$A,A'\in\cS_r$. Let $(a,b)$ and $(a',b')$ denote respectively rows
of $A$ and $\tA'$ with largest norm.  Since $a\ge rb$, we have
\begin{equation*}
 \|AA'\| \ge aa'+bb' \ge rb(a'+b') > r b\|A'\|.
\end{equation*}
Similarly, since $a'\ge rb'$, we find $\|AA'\| > r b'\|A\|$. If
$b=\|A\|$ or $b'=\|A'\|$, this gives $\|AA'\| > r\|A\|\|A'\|$ as
requested.  Otherwise, we have $a=\|A\|$ and $a'=\|A'\|$ and we get
the stronger inequality $\|AA'\| > \|A\| \|A'\|$.
\end{proof}

The next proposition is more specific than Theorem
\ref{results:thmD} and thereby proves it.

\begin{proposition}
Put $W_1=\begin{pmatrix} m &m\\ m-1 &m\end{pmatrix}$ and $W_2 =
\begin{pmatrix} 2m &m\\ 2m-1 &m\end{pmatrix}$ for a non-zero
integer $m$.  Then the Fibonacci sequence $(W_i)_{i\ge 1}$ of $\cP$
generated by these two matrices is associated to a badly
approximable real number $\xi$, and it satisfies $\det W_i = m$ for
each index $i$ which is not divisible by $3$.  If $|m|$ is not the
square of an integer, then $\xi$ is not conjugate under the action
of $\GL_2(\bQ)$ to an extremal real number having an associated
Fibonacci sequence in $\GL_2(\bZ)$.
\end{proposition}

\begin{proof}
A short computation shows that $W_1$, $W_2$, $W_1W_2$ and $W_2W_1$
are linearly independent over $\bQ$.  Then, $W_1$ and $W_2$ fulfill
the hypotheses of lemma \ref{B:lemma2} and so the sequence
$(W_i)_{i\ge 1}$ is admissible.  One can check that a corresponding
matrix $N$ is $\begin{pmatrix} m &-m\\ -2m &2m-1\end{pmatrix}$.
Moreover, for the given choice of $m$, the matrices $W_1$ and $W_2$
satisfy the hypotheses of Lemma \ref{D:lemma1} with
$W=\begin{pmatrix} 0 &0\\ -1 &0\end{pmatrix}$. Thus, defining the
map $\Phi\colon E^*\to\cP$ as in this lemma, we have $\det
\Phi(u)=1$ for each word $u\in E^*$ of even length and $\det
\Phi(u)=m$ for each $u\in E^*$ of odd length. Since the length of
$w_i$ is even if and only if $i$ is divisible by $3$, we deduce that
$W_i=\Phi(w_i)$ has determinant $1$ when $i$ is divisible by $3$ and
determinant $m$ otherwise.  In particular, we have $|\det W_i|\le
|m|$ for each $i\ge 1$.

A short computation also gives $W_3 = \pm \begin{pmatrix} 3m-1 &3m\\
3m-2 &3m-1\end{pmatrix}$ and shows, in the notation of Lemma
\ref{D:lemma2}, that $\pm W_2$ and $\pm W_3$ both belong to
$\cS_{1/2}$ for some appropriate choice of signs.  Thus, for each
$i\ge 2$, one of the matrices $\pm W_i$ belongs to $\cS_{1/2}$ and
we have
\begin{equation*}
 \|W_{i+1}W_i\| > \frac{1}{2} \|W_{i+1}\| \|W_i\|.
\end{equation*}
Since the determinant of $W_{i+1}W_i$ is a divisor of $m^2$, the
content of this product is a divisor of $m$ and so the matrix
$W_{i+2} = (W_{i+1}W_i)^\red$ satisfies $\|W_{i+2}\| \ge |m|^{-1}
\|W_{i+1}W_i\|$.  Combining this inequality with the previous one,
we deduce that
\begin{equation*}
 \|W_{i+2}\| > \frac{1}{2|m|} \|W_{i+1}\| \|W_i\|,
\end{equation*}
for each $i\ge 2$.  By induction, this implies $\|W_{i+1}\| >
\|W_i\| \ge 2 |m|$ for each $i\ge 2$, and so the sequence
$(W_i)_{i\ge 1}$ is unbounded.  Applying Theorem \ref{results:thmA},
we deduce that the sequence $(W_i)_{i\ge 1}$ is associated to some
extremal real number $\xi$.  Moreover, since we have $|\det
\Phi(u)|\le |m|$ for each $u\in E^*$, Theorem \ref{results:thmC}
shows that $\xi$ is badly approximable.

Finally, suppose that $\xi$ is $\GL_2(\bQ)$-conjugate to an extremal
real number $\eta$ with an associated Fibonacci sequence
$(W'_i)_{i\ge 1}$ in $\GL_2(\bZ)$.  Then, there exists a matrix
$A\in\cP$ such that $(\eta,-1)$ is proportional to $(\xi,-1)A$ and,
upon denoting by $B$ the inverse of $A$ in $\cP$, we find that
$(A*W_i'*B)_{i\ge 1}$ is a Fibonacci sequence in $\cP$ which is
associated to $\xi$. So, by Theorem \ref{results:thmA}, the
sequences $(W_i)_{i\ge 1}$ and $(A*W_i'*B)_{i\ge 1}$ differ only up
to their first terms and up to multiplication by a Fibonacci
sequence in $\{-1,1\}$.  Comparing determinants, this implies that
$|m|$ is the square of an integer.
\end{proof}

\begin{remark}
The Fibonacci sequence $(W_i)_{i\ge 1}$ in $\cP$ starting with
\begin{equation*}
\begin{array}{lll}
 W_1=\begin{pmatrix} 1 &0\\ 0 &2\end{pmatrix},
 \quad
 &W_2=\begin{pmatrix} 0 &1\\ 2 &0\end{pmatrix},
 \quad
 &W_3=\begin{pmatrix} 0 &1\\ 1 &0\end{pmatrix},
 \\ \\
 W_4=\begin{pmatrix} 2 &0\\ 0 &1\end{pmatrix},
 \quad
 &W_5=\begin{pmatrix} 0 &2\\ 1 &0\end{pmatrix},
 \quad
 &W_6=\begin{pmatrix} 0 &1\\ 1 &0\end{pmatrix},
\end{array}
\end{equation*}
is periodic of period $6$ as one finds that $W_7=W_1$ and $W_8=W_2$.
Therefore, $(W_i)_{i\ge 1}$ is a Fibonacci sequence of matrices with
bounded determinant.  It does not correspond to an extremal real
number as the sequence itself is bounded.  However, if $\Phi\colon
E^*\to\cP$ denotes the morphism of monoids sending $w_i$ to $W_i$
for each $i\ge 1$, then, for each $i\ge 1$, the word
$v_i=w_{6i+1}\cdots w_7 w_1$ is a prefix of $f_{a,b}$ whose image
under $\Phi$ is the matrix $W_1^{i+1}$ which has determinant
$2^{i+1}$ tending to infinity with $i$.
\end{remark}



\begin{thebibliography}{99}



\bibitem{Ca}
  J.~W.~S.~Cassels,
  {\it An introduction to Diophantine approximation},
  Cambridge U.\ Press, 1957.

\bibitem{DSa}
   H.~Davenport, W.~M.~Schmidt,
   Dirichlet's theorem on diophantine approximation,
   {\it Symposia Mathematica}, Vol. IV (INDAM, Rome, 1968/69),
   pp. 113--132;
   Academic Press, London, 1970.


\bibitem{DS}
  H.~Davenport, W.~M.~Schmidt,
  Approximation to real numbers by algebraic integers,
  {\it Acta Arith.\ }{\bf 15} (1969),
  393--416.

\bibitem{Lo}
  M.~Lothaire,
  {\it Combinatorics on words},
  Encyclopedia of mathematics and its applications, vol.~17,
  Addison-Wesley Pub.\ Co., 1983.

\bibitem{Lu}
  B.~Lucier,
  Binary morphisms to ultimately periodic words,
  manuscript,
  arXiv:0805.1373v1 [cs.DM].

\bibitem{Ra}
  D.~Roy,
  Approximation simultan\'ee d'un nombre et de son carr\'e,
  {\it C.\ R.\ Acad.\ Sci., Paris, ser.\ I} {\bf 336} (2003),
  1--6.

\bibitem{Rb}
  D.~Roy,
  Approximation to real numbers by cubic algebraic integers I,
  {\it Proc.\ London Math.\ Soc.\ }{\bf 88} (2004), 42--62.

\bibitem{Rd}
  D.~Roy,
  Diophantine approximation in small degree,
  in: {\it Number theory}, E.~Z.~Goren and H.Kisilevsky Eds,
  CRM Proceedings and Lecture Notes {\bf 36}
  (Proceedings of CNTA-7),
  2004, 269--285;
  arXiv:math.NT/0303150.

\bibitem{Re}
  D.~Roy,
  On two exponents of approximation related to a real number and
  its square, {\it Canad.\ J.\ Math.\ }{\bf 59} (2007), 211--224.

\bibitem{Sc}
  W.~M.~Schmidt,
  {\it Diophantine approximation}, Lecture Notes in Math.,
  vol.~785, Sprin\-ger-Verlag, 1980.

\end{thebibliography}
\end{document}